\documentclass[12pt]{article}

\usepackage{amsmath,amssymb,amsthm,amscd}
\usepackage[colorlinks, linkcolor=blue, anchorcolor=gree, citecolor=red]{hyperref}
\usepackage[numbers,sort&compress]{natbib}
\usepackage{graphicx}

\begin{document}

\newcommand{\Cyc}{{\rm{Cyc}}}\newcommand{\diam}{{\rm{diam}}}

\newtheorem{thm}{Theorem}[section]
\newtheorem{pro}[thm]{Proposition}
\newtheorem{lem}[thm]{Lemma}
\newtheorem{fac}[thm]{Fact}
\newtheorem{cor}[thm]{Corollary}
\theoremstyle{definition}
\newtheorem{ex}[thm]{Example}
\newtheorem{ob}[thm]{Observtion}
\newtheorem{remark}[thm]{Remark}
\newcommand{\bth}{\begin{thm}}
\renewcommand{\eth}{\end{thm}}
\newcommand{\bex}{\begin{examp}}
\newcommand{\eex}{\end{examp}}
\newcommand{\bre}{\begin{remark}}
\newcommand{\ere}{\end{remark}}

\newcommand{\bal}{\begin{aligned}}
\newcommand{\eal}{\end{aligned}}
\newcommand{\beq}{\begin{equation}}
\newcommand{\eeq}{\end{equation}}
\newcommand{\ben}{\begin{equation*}}
\newcommand{\een}{\end{equation*}}

\newcommand{\bpf}{\begin{proof}}
\newcommand{\epf}{\end{proof}}
\renewcommand{\thefootnote}{}
\newcommand{\sdim}{{\rm sdim}}

\def\beql#1{\begin{equation}\label{#1}}
\title{\Large\bf The strong metric dimension of the power graph of a finite group}

\author{{\sc Xuanlong Ma$^{1,3}$~~ Min Feng$^{2,}$\footnote{Corresponding author.}~~ Kaishun Wang$^3$
}\\[15pt]
{\small\em $^1$College of Mathematics and Information Science, Guangxi University, Nanning, 530004, China}\\
{\small\em $^2$School of Science, Nanjing University of Science and Technology, Nanjing, 210094, China}\\
{\small\em $^3$Sch. Math. Sci. {\rm \&} Lab. Math. Com. Sys., Beijing Normal University, Beijing, 100875, China}\\
}

 \date{}

\maketitle

\begin{abstract}
We characterize  the strong metric dimension of the power graph of a finite group. As applications, we compute the strong
metric dimension of the power graph of
a cyclic group, an abelian group,
a dihedral group or a generalized quaternion group.

\end{abstract}



\footnote{E-mail addresses: xuanlma@mail.bnu.edu.cn (X. Ma), fgmn\_1998@163.com (M. Feng), \\
wangks@bnu.edu.cn (K. Wang).}
\section{Introduction}
Given a graph $\Gamma$, denote by $V(\Gamma)$ and $E(\Gamma)$ the vertex set and edge set of $\Gamma$, respectively.
For $x,y,z\in V(\Gamma)$, we say that $z$ {\em strongly resolves} $x$ and $y$ if there exists a shortest path from $z$ to $x$ containing $y$, or a  shortest path from $z$ to $y$ containing $x$. A subset $S$ of $V(\Gamma)$ is a {\em strong
resolving set} of $\Gamma$ if every pair of vertices of $\Gamma$ is strongly resolved by some vertex of $S$. The smallest cardinality
of a strong resolving set of $\Gamma$ is called the {\em strong metric dimension} of $\Gamma$ and is denoted by $\sdim(\Gamma)$.

In the 1970s, metric dimension was first introduced, independently by Harary and Melter \cite{HM76} and by Slater \cite{S75}.
This parameter has appeared in various applications (see \cite{BC11} and \cite{CH07} for more information).
In 2004, Seb\H{o} and Tannier \cite{ST04} introduced
the strong metric dimension of a graph and presented some applications of strong resolving sets to combinatorial searching.
The strong metric dimension of
corona product graphs, rooted product graphs  and strong products of graphs were studied in \cite{KYR1,KYR2,KYR3}, respectively.
The problem of computing strong metric dimension is NP-hard \cite{KYR3}.
Some theoretical results, computational approaches and recent results on  strong metric dimension can be found in \cite{KKCM14}.

The {\em power graph} $\Gamma_G$ of a finite group $G$ has the vertex set $G$ and two distinct elements are adjacent if one is a power of the other. In 2000,
Kelarev and  Quinn \cite{n1} introduced
the concept of a power graph. Recently, many interesting results on power graphs have been obtained, see \cite{acam,Cam,CGh,CGS,FMW,FMW1,kel21,kel2,kel22,MF,MFW}. A detailed list of results and open questions on power graphs can be found in \cite{AKC}.

 This paper is organized as follows. In Section \ref{Sec2},
we express the strong metric dimension of a graph with  diameter two in terms of the clique number of its reduced graph. 
Sections \ref{Sec3} and \ref{Sec4} study the clique number of  the reduced graph of the power graph of a finite group $G$. Therefore, the strong metric dimension of $\Gamma_G$ is characterized.
In Section \ref{Sec5}, we compute the strong 
metric dimension of the power graph of
a cyclic group, an abelian group,
a dihedral group or a generalized quaternion group.

\section{Reduced graphs}\label{Sec2}

Let $\Gamma$ be a connected graph.
The {\em distance} $d_{\Gamma}(x,y)$ between vertices $x$ and $y$ is the length of a shortest path from $x$ to $y$ in $\Gamma$.
The {\em closed neighborhood } of $x$ in $\Gamma$, denoted by $N_\Gamma[x]$, is the set of vertices which have distance at most one from $x$. The greatest distance between any two vertices in $\Gamma$ is called the {\em diameter} of $\Gamma$. A subset of $V(\Gamma)$ is
a {\em clique} if any two distinct vertices in this subset are adjacent
in $\Gamma$. The {\em clique number} $\omega(\Gamma)$ is the maximum cardinality of a clique in $\Gamma$.

\begin{pro}\label{basic1}
  Let $\Gamma$ be a connected graph with diameter two. Then a subset $S$ of $V(\Gamma)$ is a strong resolving set of $\Gamma$ if and only if the following conditions hold:

{\rm(i)} $V(\Gamma)\setminus S$ is a clique in $\Gamma$;

{\rm(ii)} $N_{\Gamma}[u]\ne N_{\Gamma}[v]$ for any two distinct vertices
$u$ and $v$ of $V(\Gamma)\setminus S$.
\end{pro}
\proof
Assume that $S$ is a strong resolving set of $\Gamma$.
Let $u$ and $v$ be two distinct vertices of $V(\Gamma)\setminus S$.
Since $\Gamma$ has diameter two, we have $d_{\Gamma}(u,v)=1$ by \cite[Property 2]{KKCS12}.
Then (i) holds. If $N_{\Gamma}[u]= N_{\Gamma}[v]$, then
$d_{\Gamma}(u,w)=d_{\Gamma}(v,w)$ for any $w\in S$,
and so $u$ and $v$ can not be strongly resolved by any vertex in $S$,
a contradiction. Hence (ii) holds.

For the converse, it follows from (ii) that there exists a vertex $w$ in $N_{\Gamma}[u]\setminus N_{\Gamma}[v]$ or $N_{\Gamma}[v]\setminus N_{\Gamma}[u]$. Without loss of generality, let $w\in N_{\Gamma}[u]\setminus N_{\Gamma}[v]$. By (i), we have $w\in S$.
Note that $(w,u,v)$ is a shortest path. Therefore, $w$ strongly resolves
$u$ and $v$, as desired.
\qed

\medskip

For vertices $x$ and $y$ in a graph $\Gamma$,  define $x\equiv y$ if  $N_\Gamma[x]=N_\Gamma[y]$.
Observe  that $\equiv$ is an equivalence relation.
Let $U(\Gamma)$ be a complete set of distinct representative elements for this  equivalence relation.
The {\em reduced graph} $\mathcal{R}_\Gamma$  of $\Gamma$ has the vertex set $U(\Gamma)$ and two vertices are adjacent if they are adjacent in $\Gamma$.
We get the following result immediately from Proposition~\ref{basic1}.

\begin{thm}\label{mainthm1}
Let $\Gamma$ be a connected graph with diameter two. Then
$$
\sdim(\Gamma)=|V(\Gamma)|-\omega(\mathcal{R}_\Gamma).
$$
\end{thm}

\section{The clique number of $\mathcal{R}_{G}$}\label{Sec3}

In the remaining of this paper, we always use $G$ to denote a finite group. 
For simplify, denote by $\mathcal R_G$ the reduced graph  $\mathcal{R}_{\Gamma_G}$. Note that the diameter of $\Gamma_G$ is at most two.  
In order to compute $\sdim(\Gamma_G)$, we only need to study $\omega(\mathcal R_G)$ from Theorem~\ref{mainthm1}.

For a positive integer $n$, write
\begin{equation}\label{n}
   n=p_1^{r_1}p_2^{r_2}\cdots p_m^{r_m},
\end{equation}
where $p_1,p_2,\ldots,p_m$ are distinct prime numbers and $r_i\ge 1$
for $1\le i \le m$. Write
$$\sigma_n=\left\{
                                  \begin{array}{ll}
                                    1, & \hbox{if $m=1$;}\\
                                    \sum\limits_{i=1}^mr_i, & \hbox{if }m\ge 2.
                                  \end{array}
                                \right.$$
Let  $\mathbb{Z}_n$ be the cyclic group of order $n$.
\begin{thm}\label{cgroup}
$\omega(\mathcal R_{\mathbb Z_n})=\sigma_n$.
\end{thm}

In the rest of this section,  assume that $G$ is  a noncyclic group.
Denote by $\mathcal M$ the set of all maximal cyclic subgroups of $G$. Given a prime $p$,
let $\mathcal M_p$ be the set of all $p$-subgroups in $\mathcal M$.
Suppose $\mathcal M_p\ne \emptyset$. Let
\begin{equation}\label{mp}
 \mathcal M_p=\{M_1,M_2,\ldots, M_t\}.
\end{equation}
For $i\in\{1,\ldots,t\}$, write
$$
\{M_i\cap M_j:j\in\{1,\ldots,t\}\}=\{C_{i1},\ldots,C_{is_i}\}.
$$
Note that $C_{i1},\ldots,C_{is_i}$ are subgroups of $M_i$ which is a cyclic group of prime power order.
Without loss of generality, we may assume that
\begin{equation}\label{Mi}
 C_{i1}\subsetneq C_{i2}\subsetneq\cdots\subsetneq C_{is_i}=M_i.
\end{equation}
Let
$$
p^{\lambda_i}=\left\{\begin{array}{ll}
                       \max\{|M_i\cap M|: M\in\mathcal M\setminus\mathcal M_p\}, & \textup{if }\mathcal M_p\subsetneq\mathcal M; \\
                       p^{-1}, & \textup{if }\mathcal M_p=\mathcal M.
                     \end{array}\right.
$$
Define
\begin{equation}\label{ap}
\alpha_p=\max\{s_i-s_i'+\lambda_i+2:1\le i\le t\},
\end{equation}
where $s_i'=\min\{u:1\le u\le s_i, p^{\lambda_i}<|C_{iu}|\}$.
If $\mathcal M_p=\emptyset$, we define $\alpha_p=0$.

A finite group is called a {\em $CP$-group} \cite{He06} if every element of the group
has prime power order.
The set of all prime divisors of a positive integer $n$ is denoted by $\pi_n$.

\begin{thm}\label{mainthm2}
Let $G$ be a noncyclic group of order $n$.

{\rm (i)} If $G$ is a $CP$-group, then
$\omega(\mathcal R_G)=\max\{\alpha_p: p\in \pi_n\}.$

{\rm (ii)} If $G$ is not a $CP$-group, then
$$
\omega(\mathcal R_G)=\max\{
\{\alpha_p: p\in \pi_n\}\cup
\{\sigma_{|M|}+1: M\in\mathcal{M}\setminus (\bigcup_{p\in \pi_{n}}\mathcal{M}_p)\}
\}.
$$
\end{thm}

Combining Theorems \ref{mainthm1}, \ref{cgroup} and \ref{mainthm2}, we get a characterization of the strong metric dimension of $\Gamma_G$. 

\section{Proofs of  Theorems~\ref{cgroup}~and~\ref{mainthm2}}\label{Sec4}

\noindent{\em Proof of Theorem~\ref{cgroup}:}
Let $n$ be as in (\ref{n}).
Note that the power graph of a cyclic group of prime power order is complete. Then $\omega(\mathcal R_{\mathbb Z_n})=1$ if  $m=1$. In the following proof, we assume that $m\ge 2$.

Let $S$ be a clique of $\mathcal R_{\mathbb{Z}_n}$, and write
$S=\{u_1,u_2,\ldots,u_s\}$. We may assume that $m_1\le m_2\le\cdots\le m_s$, where  $m_i=|u_i|$.
Since $\{u_i,u_j\}\in E(\Gamma_{\mathbb{Z}_n})$ for $1\le i<j\le s$, one has $m_i|m_j$.
Note that $N_{\Gamma_{\mathbb{Z}_n}}[u_i]\ne N_{\Gamma_{\mathbb{Z}_n}}[u_j]$. It follows from \cite[Proposition~3.6]{FMW} that the elements with the same order have the same closed neighborhoods,
and the closed neighborhoods of the identity and any generator are equal.
Hence,
one gets $m_i\ne m_j$
and
$$
|S\cap \{x\in \mathbb{Z}_n: |x|=1 \text{ or } n\}|\le 1.
$$
Therefore, we get $s=|S|\le \sigma_n,$ and so
$$
\omega(\mathcal R_{\mathbb{Z}_n})\le\sigma_n.
$$

Let
$T=\{a_1,a_2,\cdots,a_{\sigma_n}\}$, where
$$
\begin{array}{l}
|a_1|=p_m,|a_2|=p_m^2,\ldots,|a_{r_m}|=p_m^{r_m},\\
|a_{r_m+1}|=p_{m-1}p_m^{r_m},\ldots,|a_{r_m+r_{m-1}}|=
p_{m-1}^{r_{m-1}}p_m^{r_m}, \\
  |a_{r_m+r_{m-1}+1}|=
p_{m-2}p_{m-1}^{r_{m-1}}p_m^{r_m},
\ldots,|a_{\sigma_n-1}|=p_1^{r_1-1}p_2^{r_2}\cdots p_m^{r_m}, |a_{\sigma_n}|=p_1^{r_1}p_2^{r_2}\cdots p_m^{r_m}.
\end{array}
$$
Then any pair of elements in $T$ has distinct closed neighborhoods in $\Gamma_{\mathbb{Z}_n}$.
Therefore, $T$ is a clique in $\mathcal R_{\mathbb{Z}_n}$ with size $\sigma_n$, as desired.
\qed

\medskip

In the following, suppose that $G$ is a noncyclic group. With reference to (\ref{mp}) and (\ref{Mi}),
for $1\le i\le t$, let $c_{iu}$ be a generator of $C_{iu}$ and $C_{i0}=\emptyset$, where $1\le u \le s_i$.

\begin{lem}\label{chain1}
Let $1\le u, v \le s_i$ and $c\in G$.

{\rm(i)} If $u\ne v$, then $N_{\Gamma_G}[c_{iu}]\ne N_{\Gamma_G}[c_{iv}]$;

{\rm(ii)}  Suppose that the maximal cyclic
subgroups of $G$ that contains $\langle c\rangle$ are
$p$-groups.  If $C_{i(u-1)}\subsetneq \langle c\rangle \subseteq C_{iu}$, then $N_{\Gamma_G}[c]=N_{\Gamma_G}[c_{iu}]$.
\end{lem}
\proof {\rm(i)} Without loss of generality, assume that $u<v$. Clearly, there exists a maximal cyclic subgroup $M_j\in \mathcal M_p\setminus\{M_i\}$ such that $C_{iu}=M_i\cap M_j$. Since $C_{iu}\subsetneq C_{iv}$, one has $C_{iv}\nsubseteq M_j$. Hence, any generator of $M_j$ belongs to $N_{\Gamma_G}[c_{iu}]\setminus N_{\Gamma_G}[c_{iv}]$, which implies that $(1)$ holds.

{\rm(ii)} We only need to show that $x\in N_{\Gamma_G}[c]$ is equivalent to $x\in N_{\Gamma_G}[c_{iu}]$. Note that $\langle c\rangle \subseteq C_{iu}$ and $C_{iu}$ is a cyclic $p$-group. Hence, it suffices to prove that if
$\langle c\rangle \subseteq \langle x\rangle$, then $x\in N_{\Gamma_G}[c_{iu}]$. Let $M$ be a maximal cyclic subgroup of $G$ that contains $\langle x\rangle$. Then $\langle c\rangle\subseteq M$, and so $M\in\mathcal M_p$. Write $M_j=M$.  Then $\langle c\rangle\subseteq M_i\cap M_j$. By (\ref{Mi}), one
has $C_{iu}\subseteq M_i\cap M_j$.
Hence, both $\langle x\rangle$ and $C_{iu}$ are subgroup of the cyclic $p$-group $M_j$, and so $x\in N_{\Gamma_G}[c_{iu}]$, as wanted.
\qed

\medskip

\begin{lem}\label{lap}
Let $p$ be a prime number and $\{a_1,a_2,\ldots,a_k\}$ be a subset of $G$. Then $k\le \alpha_p$ if the following conditions hold:

{\rm (i)} $\langle a_1\rangle\subsetneq \langle a_2\rangle \subsetneq \cdots \subsetneq \langle a_k\rangle$;

{\rm (ii)} $|a_k|$ is a power of $p$;

{\rm (iii)} $N_{\Gamma_G}[a_u]\ne N_{\Gamma_G}[a_v]$ for $1\le u< v \le k$.
\end{lem}
\proof
By condition~(ii), there exists $M_i\in \mathcal{M}_p$ such that
$\langle a_k\rangle \subseteq M_i$.
Set $|a_k|=p^m$.
If $m\le \lambda_i$, then $k\le m+1 \le \lambda_i+1< \alpha_p$
by condition~(i) and equation~(\ref{ap}). In the following proof, suppose $m> \lambda_i$.
Let $l=\min\{j: p^{\lambda_i}< |a_j|, 1\le j \le k\}$.
Then $|a_{l-1}|\le p^{\lambda_i}$, and so
\begin{equation}\label{lli}
l-1\le \lambda_i+1.
\end{equation}
For every $a_j\in\{a_l,a_{l+1},\ldots,a_k\}$, the maximal cyclic
subgroups of $G$ containing $\langle a_j\rangle$ are
$p$-groups. As refer to (\ref{Mi}), we get $C_{i(u-1)}\subsetneq\langle a_j\rangle\subseteq C_{iu}$,
where $u\in \{s_i',s_i'+1,\ldots,s_i\}$.
If $k-l> s_i-s_i'$, then there exists  indices $j\in\{l,l+1\ldots,k\}$ and $u\in \{s_i',s_i'+1,\ldots,s_i\}$ such that
$C_{i(u-1)}\subsetneq\langle a_j\rangle\subseteq\langle a_{j+1}\rangle \subseteq C_{iu}$, and it follows from Lemma~\ref{chain1}~(ii) that $N_{\Gamma_G}[a_j]=N_{\Gamma_G}[a_{j+1}]$, contradicting  condition~(iii). Hence, one has
\begin{equation}\label{ks}
k-l\le s_i-s_i'.
\end{equation}

Combining (\ref{lli}) and (\ref{ks}) we get $k\le s_i-s_i'+\lambda_i+2$,
and so $k\le \alpha_p$.
\qed

\medskip

\begin{lem}\label{mpne}
 If $\mathcal{M}_p\ne \emptyset$,
 then $\mathcal R_{G}$ has a clique with size $\alpha_p$.
\end{lem}
\proof With reference to (\ref{mp}) and (\ref{ap}), take $M_i \in \mathcal{M}_p$ such that $\alpha_p=s_i-s_i'+\lambda_i+2$. As refer to (\ref{Mi}), we have
$$
\langle c_{i1}\rangle\subsetneq\cdots\subsetneq\langle c_{i(s_i'-1)}\rangle\subsetneq\langle c_{is_i'}\rangle\subsetneq\langle c_{i(s_i'+1)}\rangle\subsetneq\cdots\subsetneq\langle c_{is_i}\rangle.
$$
For $j\in\{0,1,\ldots,\lambda_i\}$, let $x_j$ be an element of order $p^j$ in $\langle c_{is_i'}\rangle$.
Write
$$
S=\{x_0,x_1,x_2,\ldots,x_{\lambda_i},
c_{is_i'},c_{i(s_i'+1)},\ldots,c_{is_i}\}.
$$
Then $|S|=\alpha_p$ and $S$ is a clique of $\Gamma_G$.
Now we only need to show that for any two distinct vertices $a$ and $b$ in $S$,
\begin{equation}\label{Nne}
  N_{\Gamma_G}[a]\ne N_{\Gamma}[b].
\end{equation}

If $\{a,b\}\subseteq \{c_{is_i'},c_{i(s_i'+1)},\ldots,c_{is_i}\}$,
then (\ref{Nne}) holds by Lemma~\ref{chain1}~(i). In the following, suppose that $\{a,b\}\nsubseteq \{c_{is_i'},c_{i(s_i'+1)},\ldots,c_{is_i}\}$. Without loss of generality, assume that $a=x_u$ for $0\le u\le\lambda_i$.
Let $\langle x\rangle$ be a maximal cyclic subgroup containing $x_{\lambda_i}$ in $\mathcal{M}\setminus \mathcal{M}_p$.
If $b=x_v$ for $0\le v\le\lambda_i$, then $|x|$ is not a prime power, which implies that (\ref{Nne}) holds as $\{a,b\}\subseteq\langle x\rangle$. If $b=c_{iv}$ for $s_i'\le v\le s_i$, then $x\in N_{\Gamma_G}[a]\setminus N_{\Gamma_G}[b]$, and so (\ref{Nne}) holds.
\qed

\medskip

Now we prove Theorem~\ref{mainthm2}.

\medskip

\noindent{\em Proof of Theorem~\ref{mainthm2}:}
Let $S$ be a clique $\mathcal{R}_G$ with size $\omega(\mathcal{R}_G)$, and write
$$
S=\{x_1,x_2,\ldots,x_s\}.
$$
In view of \cite[Lemma 1]{man}, we may assume that $\langle x_1\rangle \subsetneq\langle x_2\rangle \subsetneq\cdots\subsetneq\langle x_s\rangle\subseteq K,$ where $K\in\mathcal M$.

Suppose that $K\in \mathcal{M}_q$ for some $q\in \pi_n$.
Then $|x_s|$ is a power of $q$. Furthermore, we have that $\langle x_u\rangle\subsetneq\langle x_{u+1}\rangle$ and $N_{\Gamma_G}[x_u]\neq N_{\Gamma_G}[x_v]$ for any $1\le u<v\le s$. According to Lemma~\ref{lap} we get
\begin{equation}\label{sle2}
  s\le \alpha_q\le\max\{\alpha_p:p\in\pi_n\},
\end{equation}
which implies that $\omega(\mathcal{R}_G)\le \max\{\alpha_p:p\in\pi_n\}$, and so (i) holds by
Lemma~\ref{mpne}.

In the following, we prove (ii).
If $K\in\mathcal{M}\setminus (\bigcup_{p\in \pi_{n}}\mathcal{M}_p)$, then
\begin{equation}\label{sle1}
s\le \sigma_{|K|}+1 \le \max\{\sigma_{|M|}+1: M\in\mathcal{M}\setminus (\bigcup_{p\in \pi_{n}}\mathcal{M}_p)\}.
\end{equation}
Combining (\ref{sle2}) and (\ref{sle1}) we have
$$
\omega(\mathcal{R}_G)\le l,\quad\textup{where }l=\max\{
\{\alpha_p: p\in \pi_n\}\cup
\{\sigma_{|M|}+1: M\in\mathcal{M}\setminus (\bigcup_{p\in \pi_{n}}\mathcal{M}_p)\}
\}.
$$
Now it suffices to show that $\mathcal{R}_G$ has a clique with size $l$.
If $l=\alpha_{q'}$ for some $q'\in\pi_n$, it follows from
Lemma~\ref{mpne} that $\mathcal{R}_G$ has a clique with size $l$, as wanted.
Now suppose that $l=\sigma_{|M|}+1$ for some $M\in\mathcal{M}\setminus (\bigcup_{p\in \pi_{n}}\mathcal{M}_p)$.
By the proof of Theorem~\ref{cgroup}, there exists a subset $A$ of $K$ with
$e\notin A$ and $|A|=l-1$ such that $K\setminus A$ is a clique of $\mathcal{R}_{M}$, where $e$ is the identity of $G$.
Let $T=A\cup \{e\}$.
It follows from \cite[Proposition 4]{Cam}  that $T$ is a clique of $\mathcal{R}_{G}$
with size $l$. The proof is now complete.
\qed

\medskip

\section{Examples}\label{Sec5}

The following result is immediate from Theorems~\ref{mainthm1} and \ref{mainthm2}.
\begin{cor}\label{pgroup}
  Let $G$ be a noncyclic $p$-group of order $n$. With reference to {\rm(\ref{Mi})},
  $$
  \sdim(\Gamma_G)=n-\max\{s_i:1\le i\le t\}.
  $$
\end{cor}

\bre\label{uppbound}
Take $p\in\pi_{|G|}$. With reference to (\ref{mp}), (\ref{Mi}) and (\ref{ap}),
we have
\begin{equation*}
s_i\le f_i+1\quad\textup{and}\quad \alpha_p\le \max\{f_i:1\le i\le t\}+1,
\end{equation*}
where $p^{f_i}=|M_i|$.
\ere

Let $\mathbb{Z}_p^n$ be the elementary abelian $p$-group of order $p^n$.
Then any maximal cyclic subgroup of $\mathbb{Z}_p^n$ is of order $p$.
Thus, by Corollary~\ref{pgroup} we have the following result.

\begin{ex}
$\sdim(\Gamma_{\mathbb{Z}_p^n})=p^n-2$.
\end{ex}

Let $D_{2n}$ denote the dihedral group of order $2n$.

\begin{ex}
For $n\ge 3$, we have $\sdim(\Gamma_{D_{2n}})=2n-(\sigma_n+1)$.
\end{ex}
\proof
Let $D_{2n}=\langle a,b: a^n=b^2=1, bab=a^{-1}\rangle$. Then
every maximal cyclic subgroup of $D_{2n}$ is either $\langle a\rangle$ or isomorphic to $\mathbb{Z}_2$.
If $n$ is a power of $2$, then $s_i=2$ for each $i$, which implies that $\sdim(\Gamma_{D_{2n}})=2n-2=2n-(\sigma_n+1)$ from Corollary~\ref{pgroup}. If $n$ is a power of an odd prime $q$,
then $\alpha_q=\alpha_2=2$,  the required result follows form Theorems~\ref{mainthm1} and \ref{mainthm2}~(i).
If $n$ is not a power of any prime, it follows from Theorems~\ref{mainthm1} and \ref{mainthm2}~(ii) and Remark~\ref{uppbound} that $\sdim(\Gamma_{D_{2n}})=2n-(\sigma_n+1)$, as wanted.
\qed

\medskip

The generalized quaternion group  is defined by
\begin{equation*}\label{2}
Q_{4n}=\langle x,y: x^n=y^2, x^{2n}=1, y^{-1}xy=x^{-1}\rangle,\qquad n\ge 2.
\end{equation*}

\begin{ex}
$\sdim(\Gamma_{Q_{4n}})=4n-(\sigma_{2n}+1)$.
\end{ex}
\proof
Note that every maximal cyclic subgroup of $Q_{4n}$ is
either $\langle x\rangle$ or isomorphic to $\mathbb{Z}_4$,
and the intersection of any two distinct maximal cyclic subgroups is of order $2$.
If $n$ is a power of $2$, then $s_i=2$ for all $i$, which implies that $\sdim(\Gamma_{Q_{4n}})=2n-(\sigma_{2n}+1)$ by Corollary~\ref{pgroup}. If $n$ is not a power of $2$, then by Theorems~\ref{mainthm1} and \ref{mainthm2}~(ii) and Remark~\ref{uppbound}, we get the desired result.
\qed

\medskip

Let $A$ be an abelian group of non-prime-power order.
By fundamental theorem of finite abelian groups, we may assume that
$$
A\cong \mathbb{Z}_{d_1}\times \mathbb{Z}_{d_2} \times \cdots
\times \mathbb{Z}_{d_k},
$$
where $d_i|d_{i+1}$ for $1\le i \le k-1$. Note that $d_k$ is not a power of a prime, and the order of any maximal cyclic subgroup of $A$ is a divisor of $d_k$. The following result follows from Remark~\ref{uppbound} and Theorems~\ref{mainthm1} and \ref{mainthm2}~(ii).

\begin{ex}
$\sdim(\Gamma_{A})=\prod_{i=1}^kd_i-(\sigma_{d_k}+1)$.
\end{ex}

It is clear that $\sdim(\Gamma_G)=|G|-1$ if and only if $G$ is
a cyclic group of prime power order. Finally, we classify the group $G$ with $\sdim(\Gamma_G)=|G|-2$.

\begin{cor}\label{sd=n-2}
Let $G$ be a group of order $n$. Then
$\sdim(\Gamma_G)=n-2$ if and only if $G$ is one of the following groups:

{\rm(i)} $\mathbb{Z}_{pq}$, where $p$ and $q$ are two distinct prime numbers;

{\rm(ii)} $Q_{2^m}$, where $m\ge 3$;

{\rm(iii)} A noncyclic $CP$-group such that
each two maximal cyclic subgroups have trivial intersection.
\end{cor}
\proof
Suppose that $\sdim(\Gamma_G)=n-2$. If $G$ is a cyclic group, by Theorems~\ref{mainthm1} and \ref{cgroup},
one has that $G\cong\mathbb{Z}_{pq}$, where $p$ and $q$ are two distinct prime numbers.
Now suppose that $G$ is not cyclic.
If $G$ has a nonidentity element $x$ with $N_{\Gamma_G}[x]=G$,
then $G\cong Q_{2^m}$ for some $m\ge 3$ by \cite[Proposition 4]{Cam}.
Assume that any nonidentity element $x$ of $G$ satisfies  $N_{\Gamma_G}[x]\ne G$. It follows from Theorems~\ref{mainthm1} and \ref{mainthm2}~(ii) that
$G$ is a $CP$-group. If there exist two maximal cyclic subgroups $\langle x\rangle$ and $\langle y\rangle$ in $G$ such that they has nontrivial intersection $\langle z\rangle$, then take $\overline T=\{e,z,x\}$ and so $\sdim(\Gamma_G)\le n-3$ by Peoposition~\ref{basic1}, a contradiction.
Thus, the necessity is valid.

By Theorems~\ref{mainthm1}, \ref{cgroup} and \ref{mainthm2}, it is routine to check that every group $G$ in (i), (ii) and (iii) satisfies
$\sdim(\Gamma_G)=n-2$.
\qed

\section*{Acknowledgement}
This work was carried out during Min Feng's visit to the Beijing Normal University
(Jun.--Jul. 2016). Kaishun Wang's research was supported by National Natural Science Foundation of China (11271047, 11371204) and the Fundamental Research Funds for the Central University of China.

\end{document}